\documentclass[a4paper,11pt]{amsart}

\usepackage{amssymb}
\usepackage{amsmath}

\theoremstyle{plain}
\newtheorem{theorem}{Theorem}[section]
\newtheorem{corollary}[theorem]{Corollary}
\newtheorem{conjecture}[theorem]{Conjecture}
\newtheorem{lemma}[theorem]{Lemma}
\newtheorem{proposition}[theorem]{Proposition}
\theoremstyle{remark}
\newtheorem{remark}[theorem]{Remark}
\newtheorem{example}[theorem]{Example}

\theoremstyle{definition}

\newcommand{\Z}{\mathbb{Z}}
\newcommand{\Q}{\mathbb{Q}}

\newcommand{\mt}{\overset{M}{\otimes}}
\newcommand{\bmt}{\overset{M}{\bigotimes}}

\newcommand{\nrm}{\operatorname{N}}
\newcommand{\tr}{\operatorname{Tr}}
\newcommand{\res}{\operatorname{R}}
\newcommand{\tor}{\operatorname{Tor}}
\newcommand{\cre}{\operatorname{Cor}}
\newcommand{\gal}{\operatorname{Gal}}

\newcommand{\Alb}{\operatorname{Alb}}

\newcommand{\homo}{\operatorname{Hom}}

\newcommand{\spec}{\operatorname{Spec}}

\newcommand{\coker}{\operatorname{coker}}

\newcommand{\divi}{\operatorname{div}}

\begin{document}


\title{Class field theory for a product of curves \\
 over a local field}


\author{Takao Yamazaki}
\date{\today}
\address{Mathematical Institute, Tohoku University,
  Aoba, Sendai 980-8578, Japan}
\email{ytakao@math.tohoku.ac.jp}

\begin{abstract}
We prove that the the kernel of the reciprocity map
for a product of curves over a $p$-adic field
with split semi-stable reduction
is divisible.
We also consider the $K_1$ of 
a product of curves over a number field.
\end{abstract}


\maketitle

\section{Introduction}

Let $X$ be a smooth proper geometrically connected variety 
over a field $k.$ We consider the group
\[ SK_1(X) = \coker[ \oplus_{y \in X_1} K_2 k(y) 
                 \to \oplus_{x \in X_0} k(x)^* ],
\]
where $X_i$ is the set of points in X of dimension $i$,
and the map is induced by the tame symbol.

\subsection{Class field theory}
When $k$ is a local field,
the unramified class field theory for $X$ 
is formulated in terms of the reciprocity map \cite{bloch, shuji}
\[ \rho_X: SK_1(X) \to \pi_1(X)^{ab}, \]
where $\pi_1(X)^{ab}$ is the maximal abelian quotient
of the \'etale fundamental group of $X$.
An important problem about $\rho_X$ is to determine
whether or not $\ker(\rho_X)$ is divisible.
In their pioneer work,
Bloch \cite{bloch} and Saito \cite{shuji} proved that
$\ker(\rho_X)$ is divisible when $X$ is a curve.
Jannsen and Saito \cite{js} showed that,
among other results, $\ker(\rho_X)$ is divisible
if $X$ is a surface with good reduction 
(assuming the Bloch-Kato conjecture).
On the other hand, Sato \cite{sato} constructed a $K3$ surface $X$
with bad semi-stable reduction
for which $\rho_X$ has non-divisible kernel.
See also \cite{sz1, sz2, yoshida} for other results.
In this paper, we show the following.

\begin{theorem}\label{local}
Let $k$ be a finite extension of $\Q_p$.
Let $X=C_1 \times \cdots \times C_d$
where $C_1, \cdots, C_d$ are 
smooth projective geometrically connected curves over $k$
such that $C_i(k)$ is not empty.
We assume 
the the Jacobian variety $J_i$ of $C_i$ satisfies either of 
the following conditions
for $i=1, \cdots, d-1$:

\begin{quote}
(1) $J_i$ has potentially good reduction.

\noindent
(2) The special fiber of the connected component 
of the N\'eron model of $J_i$ is
an extension of an abelian variety by a split torus.
\end{quote}

\noindent
Then $\ker(\rho_X)$ is divisible.
\end{theorem}

\begin{remark}\label{kahnremark}
(i) The assumptions (1) or (2) in Theorem \ref{local}
can be weaken to the conditions (I) or (II) 
in Proposition \ref{normmackey} of \S \ref{secondproof}.
(See also Remark \ref{counterremark}.)

\noindent
(ii) By the semi-stable reduction theorem
and the norm argument,
Theorem \ref{local} implies the following statement:
if $C_1, \cdots, C_d$ are 
smooth projective geometrically connected curves over 
a finite extension $k$ of $\mathbb{Q}_p$,
then $\ker(\rho_{C_1 \times \cdots \times C_d})$ 
is divisible by almost all prime.

\noindent
(iii) Kahn \cite{kahn} proved the injectivity of the reciprocity map 
(and the vanishing of the albanese kernel)
for a product of curves over a {\it finite} field.
(This is a special case of the main result of 
Kato and Saito \cite{ks}.)
Our proof of Theorem \ref{local} (and Theorem \ref{global} below)
is very close to Kahn's proof.
\end{remark}

\subsection{Conjectures of Bloch and Kahn}
We set
\begin{equation}\label{skonenorm}
 V(X) = \ker[ SK_1(X) \to k^* ] 
\end{equation}
where the map is induced by the norm map
$k(x)^* \to k^*$ for all $x \in X_0.$
(If $X$ admits a $k$-rational point,
we have $SK_1(X) \cong V(X) \oplus k^*$.)
Bloch conjectured (\cite{bloch} Remark 1.24)
that $V(C)$ should be torsion
when $C$ is a smooth proper curve over a global field.
Here we recall a theorem of Raskind concerning this conjecture.
(See \cite{ak} for other results.)
\begin{theorem}\label{raskind}(Raskind \cite{raskind}.)
If $C$ is a smooth proper geometrically connected curve 
over a number field $k$, then
$V(C)/V(C)_{\tor}$ is a divisible group.
(Here $V(C)_{\tor}$ is the subgroup of torsion elements in $V(C).$)
\end{theorem}
\begin{remark}
(i) 
Let $\bar{k}$ be an algebraic closure of $k$.
We use the notation $X_{k'} = X \times_{\spec k} \spec k'$
for any field extension $k'/k$.
Since $V(C_{\bar{k}})$ is uniquely divisible
(\cite{raskind} Lemma 1.1)
and since $\ker[V(C) \to V(C_{\bar{k}})]$ is torsion,
Theorem \ref{raskind} can be equivalently stated as 
each of the following statements
(cf. \cite{raskind} sublemma)
\begin{itemize}
\item $V(C) \otimes \Q/\Z = 0.$
\item $V(C) \to V(C) \otimes \Q$ is surjective.
\item $V(C) \to V(C_{\bar{k}})^{\gal(\bar{k}/k)}$ is surjective.
\end{itemize}

\noindent
(ii) According to Bloch's conjecture,
the quotient $V(C)/V(C)_{\tor}$ in Theorem \ref{raskind}
should be trivial,
but it might be very difficult to prove this,
as is pointed out in the introduction of \cite{raskind}.
\end{remark}

Now let $X=C_1 \times \cdots \times C_d$
where $C_1, \cdots, C_d$ are
smooth projective geometrically connected curves over $k$
such that $C_i(k)$ is non-empty.
The presence of $k$-rational points on $C_i$
implies that $V(C_i)$ is a direct summand of $V(X)$
(see \S \ref{productofcurves} for details).
Hence we have a (non-canonical) decomposition
\begin{equation}\label{deftilde}
 V(X) \cong (\bigoplus_{i=1}^d V(C_i) ) \oplus \tilde{V}(X). 
\end{equation}
Kahn's conjecture 
\cite{kahn} (see also Conjecture \ref{kahn})
implies that $\tilde{V}(X)$ should be trivial
when $k$ is a totally complex number field.
The following theorem gives a weak evidence to this statement.

\begin{theorem}\label{global}
Let $k$ be a totally complex number field,
and let $X=C_1 \times \cdots \times C_d$
where $C_1, \cdots, C_d$ are
smooth projective geometrically connected curves over $k$
such that $C_i(k)$ is non-empty.
We assume that the Jacobian variety of $(C_i)_{k_v}$ 
satisfies the conditions (1) or (2) of Theorem \ref{local}
for any finite place $v$ of $k$ and for $i=1, \cdots, d-1$.
(Here $k_v$ is the completion of $k$ at $v$.)
Then $\tilde{V}(X)$ is a divisible group.
\end{theorem}

\begin{remark}
(i) Similarly to Raskind's result,
$\tilde{V}(X)$ should be trivial
according to Kahn's conjecture,
but the author cannot prove it.

\noindent
(ii) 
By the semi-stable reduction theorem and the norm argument,
Theorem \ref{global} implies the following statement:
if $C_1, \cdots, C_d$ are 
smooth projective geometrically connected curves over 
a number field $k$ such that $C_i(k)$ is non-empty,
then $\tilde{V}(C_1 \times \cdots \times C_d)$ 
is the direct sum of a divisible group
and a torsion group of finite exponent.
(See Corollary \ref{lastentry}.)
\end{remark}

\subsection{Norm map on $V(X)$}
One of the key ingredients in the proof of 
Theorems \ref{local} and \ref{global}
is a description of $V(X)$ due to Raskind and Spiess \cite{rs},
which we will recall in \S \ref{raskindspiess}.
Using this description,
we are easily reduce to the following proposition.
(When $k'/k$ is a finite extension,
we have the norm map $N_k^{k'}: V(X_{k'}) \to V(X)$.)
\begin{proposition}\label{normmap}
Let $k'/k$ be a finite extension of fields.
Let $C$ be a smooth projective geometrically connected curve 
over $k$ such that $C(k)$ is non-empty,
and let $J$ be the Jacobian variety of $C.$

\noindent
(i) Assume that $k$ is a finite extension of $\Q_p$,
and that $J$ satisfies 
either of the conditions (1) or (2) in Theorem \ref{local}.
Then the norm map $N_k^{k'}: V(C_{k'}) \to V(C)$ is surjective.

\noindent
(ii) Assume that $k$ is a totally complex number field,
and that $J_{k_v}$ satisfies either of 
the conditions (1) or (2) in Theorem \ref{local}
for any finite place $v$ of $k$.
Then the norm map $N_k^{k'}: V(C_{k'}) \to V(C)$ is surjective.
\end{proposition}

In \S \ref{localsection} we shall give 
two proofs of Proposition \ref{normmap} (i).
The first proof uses the class field theory for curves
due to Bloch \cite{bloch} and Saito \cite{shuji}.
The second proof is longer, but
it shows a slightly stronger result.
(See Proposition \ref{normmackey}.
This proposition is formulated
for a more general types of semi-abelian variety $G$.
Proposition \ref{normmap} (i) follows as the case $G=J$.)
This proof is based on the classical results such as 
the uniformization of an abelian variety,
Mazur's theorem on the norm map of an abelian variety,
the local class field theory and the Tate duality.
In \S \ref{counterexample}, we discuss an example of $C$
which shows the assumption of 
Proposition \ref{normmap} (i) is essential.

In \S \ref{globalsection},
we deduce Proposition \ref{normmap} (ii) from (i)
by showing an analogue of the `Hasse norm theorem' for $V(C)$.
(See Proposition \ref{hasse} for details).
Here the key ingredient is a result from
Kato's Hasse principle \cite{kato}.

\subsection{Conventions}
Throughout this paper, $k$ is a field of characteristic zero.
We fix an algebraic closure $\bar{k}$ of $k$,
and all algebraic extension of $k$ 
is supposed to be a subfield of $\bar{k}.$
We write $G_k$ for the absolute Galois group 
$\gal(\bar{k}/k)$ of $k$.

Let $A$ be an abelian group.
For a non-zero integer $n$,
we write $A[n]$ and $A/n$ for the kernel and cokernel
of the map $n: A \to A.$
We define $A_{\divi}$
to be the maximal divisible subgroup in $A.$
When a group $G$ acts on $A$,
we write $A^{G}$ and $A_{G}$ for the invariants and coinvariants
of $A$ by $G$.

\section{Review of the results of Raskind and Spiess}\label{raskindspiess}
We recall some necessary results of Somekawa \cite{somekawa}
and Raskind-Spiess \cite{rs}.
Let $G$ be a semi-abelian variety over $k$.
For any tower of finite extensions $k_2/k_1/k$,
we have the restriction map 
$\res_{k_1}^{k_2}:G(k_1) \to G(k_2)$
and the norm map 
$\nrm_{k_1}^{k_2}:G(k_2) \to G(k_1)$.
The composition $\nrm_{k_1}^{k_2} \circ \res_{k_1}^{k_2}$
is the multiplication by $[k_2:k_1].$

\subsection{Mackey tensor product and Milnor $K$-group}\label{mackey}
Let $G_1, \cdots, G_r$ be semi-abelian varieties over $k.$
Let
\[ X(G_1, \cdots, G_r) 
   = \bigoplus_{k'/k} G_1(k') \otimes \cdots \otimes G_r(k'), \]
where $k'$ runs through all tower of finite extensions of $k$.
Let $R(G_1, \cdots, G_r)$ 
be the subgroup of $X(G_1, \cdots, G_r)$ 
generated by the elements of the form
\[  \res_{k_1}^{k_2}(a_1) \otimes \cdots \otimes a_{i_0} \otimes
    \cdots \otimes \res_{k_1}^{k_2}(a_1)
~-~ a_1 \otimes \cdots \otimes \nrm_{k_1}^{k_2}(a_{i_0}) \otimes
    \cdots \otimes a_r
\]
where $k_2/k_1/k$ is a tower of finite extensions,
$1 \leq i_0 \leq d,$
$a_i \in G_i(k_1) ~(i \not= i_0),$
and $a_{i_0} \in G_{i_0}(k_2)$.
Following Raskind and Spiess \cite{rs} (3.2) 
(see also \cite{kahn2} \S 5), 
we define
\[ G_1 \mt \cdots \mt G_r (k)  = X(G_1, \cdots, G_r)/R(G_1, \cdots, G_r) \]
to be the Mackey tensor product of $G_1, \cdots, G_r.$
The class of $a_1 \otimes \cdots \otimes a_r \in \otimes_i G_i(k')$
is denoted by $(a_1, \cdots, a_r)_{k'/k}.$
The Mackey tensor product $\mt_i G_i$
have the following functorial properties
for any tower of finite extensions $k_2/k_1/k$
(cf. loc. cit. and \cite{somekawa} \S 1.3)
\begin{align*}
\res_{k_1}^{k_2}: \bmt_i G_i (k_1) \to \bmt_i G_i (k_2),
\qquad
\nrm_{k_1}^{k_2}: \bmt_i G_i (k_2) \to \bmt_i G_i (k_1).
\end{align*}
The composition $\nrm_{k_1}^{k_2} \circ \res_{k_1}^{k_2}$
is the multiplication by $[k_2:k_1].$
Note that the norm map can be written simply as
\[
\nrm_{k_1}^{k_2}(a_1, \cdots, a_r)_{k_2'/k_2} 
=(a_1, \cdots, a_r)_{k_2'/k_1}. 
\]

Somekawa \cite{somekawa}
introduced the Milnor $K$-group $K(k; G_1, \cdots, G_r)$
attached to $G_1, \cdots, G_r$.
This is defined as a quotient of $\mt G_i (k)$
by new relations coming from the Weil reciprocity.
(We do not need the explicit description of these relations.
See \cite{somekawa} for details.)
This group again has the functional properties
corresponding to $\res_{k_1}^{k_2}$ and $\nrm_{k_1}^{k_2}$.
The following example would be helpful 
to understand $K(k; G_1, \cdots, G_r)$.

\begin{example}\label{milnorK}
Let $G_1 = \cdots = G_r = \mathbb{G}_m.$
Somekawa (\cite{somekawa} Theorem 1.4) proved that
$K(k; \mathbb{G}_m, \cdots, \mathbb{G}_m)$ is 
canonically isomorphic to the usual Milnor $K$-group $K_r^M k.$
Hence we have a surjection
\[ \bmt_{1 \leq i \leq r} \mathbb{G}_m (k) \to K_r^M k, \]
which is given by 
$(a_1, \cdots, a_r)_{k'/k} \mapsto \nrm_k^{k'}( \{a_1, \cdots, a_r\} ).$
(Here $\nrm_k^{k'}$ is the norm map on the Milnor $K$-groups.
The projection formula for the Milnor $K$-groups shows 
the well-definedness of this map.)
\end{example}

\subsection{Product of curves}\label{productofcurves}
The following theorem plays an essential role in this paper.

\begin{theorem}\label{rs}
Let $C_1, \cdots, C_d$ be 
smooth projective geometrically connected curves over a field $k$
such that $C_i(k) \not= \phi$ for each $i$.
Let $J_i$ be the Jacobian variety of $C_i$,
and let $X = C_1 \times \cdots \times C_d.$
Then there is an isomorphism 
(depending on the choice of a $k$-rational point on $C_i$)
\[ SK_1(X) \cong CH^{d+1}(X, 1) \cong 
   \bigoplus_{r=0}^{d} 
   \bigoplus_{1 \leq i_1 < \cdots < i_r \leq d}
   K(k; J_{i_1}, \cdots, J_{i_r}, \mathbb{G}_m).
\]
\end{theorem}
\begin{proof}
The first isomorphism is given by \cite{landsburg} Lemma 2.8.
The second isomorphism is an easy consequence of 
Raskind-Spiess \cite{rs}.
Here we include a brief account.
We use the notations in loc. cit.
By the decomposition of the Chow motive of $C_i$
explained in loc. cit. \S 2.4,
we have a decomposition
\[ CH^{d+1}(X, 1) \cong 
   \bigoplus_{r=0}^{d} 
   \bigoplus_{1 \leq i_1 < \cdots < i_r \leq d}
   CH^{r+1}(C_{i_1}^+ \otimes \dots \otimes C_{i_r}^+, 1).
\]
Then loc. cit. Remark 2.4.2 (c)
shows an isomorphism
\[ CH^{r+1}(C_{i_1}^+ \otimes \dots \otimes C_{i_r}^+, 1)
   \cong K(k; \mathcal{CH}^1(C_{i_1}), \cdots, \mathcal{CH}^1(C_{i_r}), 
    \mathcal{CH}^1(\spec k, 1) ),
\]
and the last group is isomorphic to
$K(k; J_{i_1}, \cdots, J_{i_r}, \mathbb{G}_m)$
by the same argument as loc. cit. \S 2.4.
\end{proof}

The projection 
$SK_1(X) \to K(k; \mathbb{G}_m) = k^*$
(the term with $r=0$)
is nothing other than the map appearing in \eqref{skonenorm}.
Applying the theorem to the case $d=1,$ 
we have
\begin{equation}
\label{surj1}
J_i \mt \mathbb{G}_m (k) \twoheadrightarrow 
K(k; J_i, \mathbb{G}_m) \cong V(C_i).
\end{equation}
We then obtain
\begin{equation}
\label{surj2}
\underset{r \geq 2; ~1 \leq i_1 < \cdots < i_r \leq d}{\bigoplus}
   J_{i_1} \mt \cdots \mt J_{i_r} 
            \mt \mathbb{G}_m (k)
 \twoheadrightarrow 
\underset{r \geq 2; ~1 \leq i_1 < \cdots < i_r \leq d}{\bigoplus}
   K(k; J_{i_1}, \cdots J_{i_r}, \mathbb{G}_m)
   \cong \tilde{V}(X)
\end{equation}
This explains the decomposition \eqref{deftilde}.

Here we recall the conjecture of Kahn.
\begin{conjecture}\label{kahn}(Kahn \cite{kahn}.)
Assume $k$ is finitely generated over its prime subfield $k_0$.
According as the characteristic of $k$ is zero or positive,
we set $s= \tr \deg(k/k_0)+1$ or $s= \tr \deg(k/k_0)$.
Let $G_1, \cdots, G_r$ be semi-abelian varieties over $k$.
If $r = s+1$, then $K(k; G_1, \cdots, G_r)$ is torsion.
If $r > s+1$ and if $k$ does not admit an ordered field structure,
then $K(k; G_1, \cdots, G_r)$ is trivial.
\end{conjecture}

Kahn \cite{kahn} proved the conjecture when $k$ is a finite field,
and applied it to the class field theory for a product of curves.
(See Remark \ref{kahnremark} (iii).
The role of Proposition \ref{normmap} is
replaced by Lang's theorem \cite{lang} there.)
In view of eq. \eqref{surj1},
Conjecture \ref{kahn} for $K(k; J, \mathbb{G}_m)$
with $k$ a global field and $J$ the Jacobian variety of a curve
is the same as Bloch's conjecture recalled in the introduction.
Similarly, with the notation in Theorem \ref{rs},
Conjecture \ref{kahn} implies $\tilde{V}(X)=0$
when $k$ is a totally complex number field.

\subsection{Surjectivity of the norm map and divisibility}
Here is a simple lemma.

\begin{lemma}
Let $V$ be a Mackey tensor product of
(any number of) semi-abelian varieties over a field $k$.
We assume the following condition
\[(*)_V \quad \nrm_{k_1}^{k_2}: V(k_2) \to V(k_1)~
\text{is surjective for any tower of finite extensions} ~k_2/k_1/k.
\]
Then $G \mt V (k)$ is a divisible group
for any semi-abelian variety $G$ over $k$.
In particular, $V'= G \mt V$ satisfies the condition $(*)_{V'}$.
\end{lemma}
\begin{proof}
Take any finite extension $k_1/k$
and any elements $a \in G(k_1)$ and $b \in V(k_1)$.
We show that $(a, b)_{k_1/k} \in G \mt V(k)$ 
is divisible by any natural number $n$.
There exists a finite extension $k_2/k_1$ and
$a' \in G(k_2)$ such that $n a' = \res_{k_1}^{k_2}(a).$
By assumption, there also exists $b' \in V(k_2)$ 
such that $\nrm_{k_1}^{k_2}(b') = b.$
Then $(a, b)_{k_1/k} = n (a', b')_{k_2/k}$,
and we are done.
\end{proof}

\begin{corollary}\label{divlemma}
Let $G_1, \cdots, G_r$ be semi-abelian varieties over $k$.
If $\nrm_{k_1}^{k_2}: 
K(k_2; G_r, \mathbb{G}_m) \to K(k_1; G_r, \mathbb{G}_m)$
is surjective for any tower of finite extensions $k_2/k_1/k$,
then 
$K(k; G_1, \cdots, G_r, \mathbb{G}_m)$
is divisible.
\end{corollary}

By this corollary and eq. \eqref{surj1} \eqref{surj2},
the proof of Theorem \ref{global}
is reduced to Proposition \ref{normmap} (ii).
As for Theorem \ref{local},
we follow the argument of Kahn \cite{kahn}.

\vspace{3mm}
\noindent
{\it Reduction of Theorem \ref{local} to Proposition \ref{normmap} (i).}
(Cf. Kahn \cite{kahn} p.1041.)
We write $\pi_1(X)^{ab, geo}$ for the kernel of
the map $\pi_1(X)^{ab} \to \pi_1(\spec k)^{ab} \cong G_k^{ab}$
induced by the structure morphism $X \to \spec k.$
We have a commutative diagram
\[
\begin{matrix}
&SK_1(X)& \cong 
  &k^*& \oplus &(\bigoplus_{i=1}^d V(C_i) )& \oplus &\tilde{V}(X)&
\\
& {}^{\rho_X}\downarrow&  &{}^{\rho_{\spec k}}\downarrow&
& {}^{\oplus\rho_{C_i}}\downarrow&  & &
\\
&\pi_1(X)^{ab}& \cong 
     &\pi_1(\spec k)^{ab}& 
     \oplus & \bigoplus_{i=1}^d \pi_1(C_i)^{ab, geo}&
\end{matrix}
\]
Here the upper row is induced by Theorem \ref{rs}
(and comments after that).
The lower row is given by
$\pi_1(X)^{ab} \cong 
 \pi_1(\spec k)^{ab} \oplus \pi_1(X)^{ab, geo}$
(deduced by the presence of a $k$-rational point),
and by
$\pi_1(X)^{ab, geo} \cong (T \Alb_X)_{G_k}
  \cong \oplus_i (T J_i)_{G_k}
  \cong \oplus_i \pi_1(C_i)^{ab, geo}.$
(Here, for an abelian variety $A$, 
we denote by $T A$ for the full Tate module of $A$.)
The map $\rho_{\spec k}$ is injective by
the usual local class field theory.
The kernel of $\rho_{C_i}$ is divisible
by the main result of the class field theory for curves \cite{shuji}.
By Proposition \ref{normmap} (i) and Corollary \ref{divlemma},
$\tilde{V}(X)$ is divisible.
This finishes the proof.
\qed

\section{Local field}\label{localsection}
In this section, 
$k$ is a finite extension of $\mathbb{Q}_p.$
We give two proofs of Proposition \ref{normmap} (i).

\subsection{First proof of Proposition \ref{normmap} (i)}\label{firstproof}
By the norm argument, we see that
$V(C_{k'})_{\divi}$ surjects onto $V(C)_{\divi}$
by $N_k^{k'}$.
Hence it is enough to show the surjectivity of the map 
$V(C_{k'})/V(C_{k'})_{\divi} \to V(C)/V(C)_{\divi}$
induced by $N_k^{k'}$.
The main result of the class field theory for curves \cite{shuji}
shows that there is a commutative diagram with exact rows
\[
\begin{matrix}
0 \to & V(C_{k'})/V(C_{k'})_{\divi}& \to &TJ_{G_{k'}}& \to 
&\hat{\Z}^{\oplus r}& \to 0 \\
     & \downarrow&  & \downarrow & &\downarrow \\
0 \to & V(C)/(C)_{\divi}& \to  & TJ_{G_{k}}& \to 
& \hat{\Z}^{\oplus r}& \to 0.
\end{matrix}
\]
Here, $r$ is the rank of $C$ 
defined in \cite{shuji} Definition 2.5.
The crucial point is that, under our assumption,
the rank of $C_{k'}$ is the same as that of $C$
by \cite{shuji} Theorem 6.2.
(We remark also that $r=0$ 
when $J$ has potentially good reduction.)
The middle vertical map is the projection
(induced by the identity map on $TJ$), hence surjective.
This shows the surjectivity of the right vertical map,
hence it is injective as well.
This shows the surjectivity of the left vertical map,
and we are done.

\subsection{Second proof of Proposition \ref{normmap} (i)}
\label{secondproof}
The following proposition implies Proposition \ref{normmap} (i),
since the theory of the rigid uniformization 
(\cite{bx} Theorem 1.2) shows that
the conditions (1) and (2) of Theorem \ref{local} are covered 
by the conditions (I) and (II) below.

\begin{proposition}\label{normmackey}
Let $G$ be a semi-abelian variety 
over a finite extension $k$ of $\Q_p$.
Assume that $G$ satisfies either of the following conditions.

\renewcommand{\labelenumi}{(\Roman{enumi})}
\begin{enumerate}
\item 
$G$ is an extension of 
an abelian variety with potentially good reduction 
by a split torus.
\item 
The $G_k$-module $G(\bar{k})$ fits in the exact sequence
$0 \to M \to E(\bar{k}) \to G(\bar{k}) \to 0$
where $E$ is a semi-abelian variety satisfying the condition (I) above
and $M$ is a split lattice.
\item 
$G$ is a torus.
\end{enumerate}
\renewcommand{\labelenumi}{(\arabic{enumi})}

Then, the norm map
$\nrm_k^{k'}: G \mt \mathbb{G}_m (k') \to G \mt \mathbb{G}_m (k)$
is surjective for any finite extension $k'/k.$
\end{proposition}

\begin{proof}
We may assume $k'/k$ is a Galois extension 
whose degree is a prime number $l$.
By the norm argument,
we may assume that $k$ contains all the $l$-th roots of unity.
We take a finite extension $k_1/k$
and elements $a \in G(k_1), b \in k_1^*.$
We shall show that $(a, b)_{k_1/k}$ is in the image of $\nrm_{k}^{k'}.$
If $k_1$ contains $k',$ 
then $(a, b)_{k_1/k} = \nrm_{k}^{k'}( (a, b)_{k_1/k'}).$
Hence we may assume  
$k_1'= k_1 \cdot k'/k_1$ is an extension of degree $l$.
The rest of the proof consists of six steps.

\vspace{2mm}
\noindent
{\it Step 1: The case where
$G$ is an abelian variety with good reduction.}
If $k_1'/k_1$ is an unramified extension,
Mazur's result (\cite{mazur} Corollary 4.4) shows that
there exists an $a' \in G(k'_1)$ such that
$N_{k_1}^{k'_1}(a')=a,$
and we have 
$(a, b)_{k_1/k} = \nrm_{k}^{k'}( (a', \res_{k_1}^{k'_1}(b))_{k_1'/k'}).$
Assume $k_1'/k_1$ is a ramified extension.
Let $k_2$ be the unramified extension of $k_1$ of degree $l$,
and let $k'_2 = k_2 \cdot k'.$
Again by Mazur's result, there exists an $a' \in G(k_2)$ such that
$N_{k_1}^{k_2}(a')=a.$
On the other hand, the map 
$k^*_2/\nrm_{k_2}^{k'_2}(k'^*_2) \to k^*_1/\nrm_{k_1}^{k'_1}(k'^*_1)$
induced by the norm map $\nrm_{k_1}^{k_2}$
is, by the local class field theory,
isomorphic to the restriction map
$\gal(k'_2/k_2) \to \gal(k'_1/k_1),$
which is
an isomorphism between cyclic groups of degree $l$
by assumption.
Since $\nrm_{k_1}^{k_2} \circ R_{k_1}^{k_2} = l$,
the restriction map $R_{k_1}^{k_2}$
induces the zero map
$k^*_1/\nrm_{k_1}^{k'_1}(k'^*_1) \to k^*_2/\nrm_{k_2}^{k'_2}(k'^*_2).$
This means the existence of $b' \in k'^*_2$ such that
$\nrm_{k_2}^{k'_2}(b') = \res_{k_1}^{k_2}(b).$
Now we have 
$(a, b)_{k_1/k} 
= \nrm_{k}^{k'}( (\res_{k_2}^{k'_2}((a'), b')_{k_2'/k'}).$

\vspace{2mm}
\noindent
{\it Step 2: The case where $G$ is a split torus.}
We may assume $G=\mathbb{G}_m.$
We take an extension $k_2/k_1$ of degree $l$,
which is linearly disjoint from $k_1.$
(The existence of such an extension follows
from the assumption $\mu_l \subset k$ 
and the Kummer theory).
Let $k'_2 = k' \cdot k_2.$
By the same argument as Step 1,
we see that
$\nrm_{k_1}^{k_2}$ induces an isomorphism 
$k^*_2/\nrm_{k_2}^{k'_2}(k'^*_2) \to k^*_1/\nrm_{k_1}^{k'_1}(k'^*_1)$
and that 
$\res_{k_1}^{k_2}$ induces the zero map 
$k^*_1/\nrm_{k_1}^{k'_1}(k'^*_1) \to k^*_2/\nrm_{k_2}^{k'_2}(k'^*_2)$.
Thus we have $a_1 \in k'^*_1, ~a_2 \in k_2^*, ~b' \in k'^*_2$
such that $a = \nrm_{k_1}^{k'_1}(a_1) \nrm_{k_1}^{k_2}(a_2),~
\res_{k_1}^{k_2}(b) = \nrm_{k_2}^{k'_2}(b').$
Now we have 
$(a, b)_{k_1/k} 
= \nrm_{k}^{k'}( (a_1, \res_{k_1}^{k'_1}(b))_{k_1'/k'}
                +(\res_{k_2}^{k'_2}(a_2), b')_{k_2'/k'}).$

\vspace{2mm}
\noindent
{\it Step 3: The case (III).}
We take a finite Galois extension $K/k$
such that $G$ becomes a split torus over $K.$
We proceed by the induction on $[K:k].$
Let $K_1 = k_1 \cdot K.$

If $K_1/k_1$ is linearly disjoint from $k_1'/k_1$,
by the argument same as Step 2 we have 
$b = \nrm_{k_1}^{k'_1}(b_1) \nrm_{k_1}^{K_1}(b_2)$
for some $b_1 \in k'^*_1, b_2 \in K_1^*$.
By Step 2, there exists $x \in G \mt \mathbb{G}_m (K'_1)$ such that 
$(\res_{k_1}^{K_1}(a), b_2)_{K_1/K_1} = \nrm_{K_1}^{K'_1}(x),$
where $K'_1 = k'_1 \cdot K_1.$
Then we have
$(a, b)_{k_1/k} 
= \nrm_{k}^{k'}( (\res_{k_1}^{k'_1}(a_1), b_1)_{k_1'/k'}
                +\nrm_{k'}^{K'_1}(x) ).$

We assume $k_1' \subset K_1.$
Let $M = \homo(G, \mathbb{G}_m)$ be the character group of $G$,
on which $G_{k}$ acts through $\gal(K/k).$
Since $H^1(K_1, M)=0$,
we have an exact sequence
\begin{equation}\label{tatedual}
 0 \to H^2(K_1/k_1, M) \to H^2(k_1, M) \to H^2(K_1, M). 
\end{equation}
We define $N$ to be the order of the finite group $H^2(K_1/k_1, M)$.
(Note that this group is dual to $G(k_1)/ \nrm_{k_1}^{K_1} G(K_1)$
under the Tate duality.)
We take an extension $k_2/k_1$ of degree $N$
which is linearly disjoint from $k_1'/k_1.$
Let $K_2 = K_1 \cdot k_2$ and $k_2' = k' \cdot k_2.$

If $[K_2: k_2] < [K_1 : k_1] (\leq [K:k])$,
the norm map 
$\nrm_{k_2}^{k_2'}: 
G \mt \mathbb{G}_m(k_2') \to G \mt \mathbb{G}_m(k_2)$
is surjective by the inductive hypothesis.
We can write 
$b = \nrm_{k_1}^{k'_1}(b_1) \nrm_{k_1}^{k_2}(b_2),~
(\res_{k_1}^{k_2}(a), b_2)_{k_2/k_2} = \nrm_{k_2}^{k_2'}(x)$
with some 
$b_1 \in k'^*_1, b_2 \in k_2^*, x \in G \mt \mathbb{G}_m(k_2')$.
Then we have
$(a, b)_{k_1/k} 
= \nrm_{k}^{k'}( (\res_{k_1}^{k'_1}(a), b_1)_{k_1'/k'}
+\nrm_{k'}^{k'_2}(x) ).$

Lastly, we assume $[K_2: k_2] = [K_1 : k_1].$
Then, in the commutative diagram
\[
\begin{matrix}
0 \to &H^2(K_1/k_1, M)& \to &H^2(k_1, M) \to &H^2(K_1, M)& \\
      &\downarrow^{\cong} &  &\downarrow^{\res^{k_2}_{k_1}}  
  &\downarrow^{\res^{K_2}_{K_1}} & \\
0 \to &H^2(K_2/k_2, M)& \to &H^2(k_2, M) \to &H^2(K_2, M)& \\
      &\downarrow^{0} &  &\downarrow^{\cre^{k_2}_{k_1}}  
  &\downarrow^{\cre^{K_2}_{K_1}} & \\
0 \to &H^2(K_1/k_1, M)& \to &H^2(k_1, M) \to &H^2(K_1, M),& 
\end{matrix}
\]
the left upper vertical map is an isomorphism
since $\gal(K_2/k_2) \cong \gal(K_1/k_1)$.
By the definition of $N,$
the lower left vertical map is the zero map.
By the Tate duality, this means that 
$R_{k_1}^{k_2}: G(k_1) \to G(k_2)$
induces the zero map
$G(k_1)/ \nrm_{k_1}^{K_1} G(K_1) \to 
   G(k_2)/ \nrm_{k_2}^{K_2} G(K_2).$
Hence there exists 
$a' \in G(K_2)$
such that 
$\res_{k_1}^{k_2}(a) = \nrm_{k_2}^{K_2}(a').$
There also exist
$b_1 \in k'^*_1$ and $b_2 \in k^*_2$
such that
$b =\nrm_{k_1}^{k'_1}(b_1)\nrm_{k_1}^{k_2}(b_2).$
Then we have
$(a, b)_{k_1/k} 
= \nrm_{k}^{k'}( (\res_{k_1}^{k'_1}(a), b_1)_{k_1'/k'}
 +(a', \res_{k_2}^{K_2}(b_2))_{K_2/k'} ).$
This completes the proof of the Case (3).

\vspace{2mm}
\noindent
{\it Step 4: The case where $G$ has potentially good reduction.}
This step is parallel to Step 3.
We only indicate how to modify the argument of Step 3.
We take $K$ to be a finite Galois extension of $k$
such that $A$ has good reduction over $K.$
In the place of eq. \eqref{tatedual},
we use an exact sequence
\[ 0 \to H^1(K_1/k_1, G^{t})
     \to H^1(k_1, G^{t})
     \to H^1(k_1, G^{t}),
\]
where $G^t$ is the dual abelian variety of $G.$
In this case, the (finite) group
$H^1(K_1/k_1, G^{t})$ is dual to 
$G(k_1)/ \nrm_{k_1}^{K_1} G(K_1)$ under the Tate duality.
The rest is straightforward.

\vspace{2mm}
\noindent
{\it Step 5: The case (I).}
We have an exact sequence $0 \to T \to G \to A \to 0$
where $T$ is a split torus and 
$A$ is an abelian variety over $k$
with potentially good reduction.
By Steps 2 and 4, the proof is reduced to showing
the exactness of 
\[ T \mt \mathbb{G}_m (k) \to G \mt \mathbb{G}_m (k)
  \to A \mt \mathbb{G}_m (k) \to 0
\]
(and a similar exactness for $k'/k$).
By Hilbert 90, the sequence
$0 \to T(k_1) \to G(k_1) \to A(k_1) \to 0$ is exact
for any $k_1/k.$
We then have the exactness of the lower row in 
the following commutative diagram
(recall the notations in \S \ref{mackey})
\[
\begin{matrix}
& & 
&R(G, \mathbb{G}_m)& \to 
&R(A, \mathbb{G}_m)&
\\
& &
&\cap&
&\cap&
\\
&X(T, \mathbb{G}_m)& \to 
&X(G, \mathbb{G}_m)& \to 
&X(A, \mathbb{G}_m)& \to 0.
\end{matrix}
\]
It remains to show the surjectivity of the upper horizontal map,
but it is immediate from the surjectivity of 
$G(k_1) \to A(k_1).$

\vspace{2mm}
\noindent
{\it Step 6: The case (II).}
We have an exact sequence 
$0 \to M \to E(\bar{k}) \to G(\bar{k}) \to 0$ of $G_k$-modules,
where $E$ is a semi-abelian variety satisfying the condition (I)
and $M$ is a free $\Z$-module of finite rank with
trivial $G_k$-action.
Since we have $H^1(k_1, M)=0$ for any finite extension $k_1/k$,
we see that $E \mt \mathbb{G}_m (k)$ 
surjects onto $G \mt \mathbb{G}_m (k).$
Hence this case follows from Step 5.
\end{proof}

\subsection{Counter example to the surjectivity}\label{counterexample}

We recall some facts from \cite{silverman} Chapter V \S 5.
Let $E$ be an elliptic curve over $k$.
(By abuse of language) we call $E$ a {\it Tate curve} 
if it has split multiplicative reduction.
If $E$ is an elliptic curve over $k$
whose $j$-invariant is not integral,
then there exists an extension $k'/k$ of degree $\leq 2$
such that $E_{k'}$ is a Tate curve over $k'$.
Here $k'/k$ is uniquely determined by $E$
unless $E$ is already a Tate curve over $k$.
We call $E$ a {\it non-split Tate curve} if $k' \not= k$.
A non-split Tate curve has
either non-split multiplicative or additive reduction over $k$.

\begin{lemma}\label{counterlemma}
Let $E$ be a non-split Tate curve 
over a finite extension $k$ of $\Q_p$.
Let $k'/k$ be a quadratic extension
such that $E_{k'}$ is a Tate curve over $k'.$
Then, $\coker[\nrm_k^{k'}: V(E_{k'}) \to V(E)]$
is isomorphic to $\Z/2\Z.$
\end{lemma}

\begin{proof}
The assumption on the reduction shows that
the rank (in the sense of \cite{shuji} Definition 2.5)
of $E$ and $E_{k'}$ is $0$ and $1$ respectively.
By the class field theory for curves (\cite{shuji}, see also \S 3.1)
we have the commutative diagram with exact upper row
\[
\begin{matrix}
0 \to & V(E_{k'})/V(E_{k'})_{\divi}& \to &TE_{G_{k'}}& \to 
&\hat{\Z}& \to 0 \\
     & \downarrow&  & \downarrow &  \\
 & V(E)/(E)_{\divi}& \cong  & TE_{G_{k}}.&  & &
\end{matrix}
\]
The right vertical arrow is surjective (cf. \S \ref{firstproof}).
Hence we have 
$\coker[\nrm_k^{k'}: V(E_{k'}) \to V(E_{k})] \cong \hat{\Z}_{\gal(k'/k)}$.
By the description of $G_k$-module structure of $E(\bar{k})$
(cf. \cite{silverman} Lemma 5.2),
the action of $\gal(k'/k)$ on $\hat{\Z}$
is seen to be non-trivial.
This completes the proof.
\end{proof}

\begin{remark}\label{counterremark}
(i)
This example demonstrates that
our method is not sufficient to prove
Theorem \ref{local} without the assumption (1) or (2).
(For example, the author does not know if
$\ker(\rho_X)$ is divisible for $X=E \times E$.)

\noindent
(ii) 
The argument of \S \ref{firstproof} shows that,
if $E_{k_1}$ is a non-split Tate curve
over a finite extension $k_1$ of $k$,
then $\nrm_{k}^{k_1}: V(E_{k_1}) \to V(E)$ is surjective.

\end{remark}

\section{number field}\label{globalsection}
Throughout this section, $k$ is a finite extension of $\Q.$

\subsection{Hasse norm theorem}
We complete the proof of Proposition \ref{normmap} (ii).
This is a consequence of Proposition \ref{normmap} (i) and
the following `Hasse norm theorem' for $V(C)$.

\begin{proposition}\label{hasse}
Let $C$ be a smooth projective geometrically connected curve 
over a number field $k$ such that $C(k)$ is not empty.
Let $k'/k$ be a finite Galois extension with the group $G$.
For each place $v$ of $k,$
we choose a place $w(v)$ of $k'$ above $v$.
Then the natural map
\[ (**) \qquad
  V(C)/\nrm_k^{k'} V(C_{k'}) \to 
   \oplus_v V(C_{k_v}) /\nrm_{k_v}^{k'_{w(v)}} V(C_{k'_{w(v)}})
\]
(here $v$ runs all places of $k$ including infinite places)
is an isomorphism.
\end{proposition}

\begin{remark}\label{satoremark}
(i) By Proposition \ref{normmap} (i),
the direct summand in the right hand side of $(**)$
is trivial for almost all $v.$

\noindent
(ii) Since $k'/k$ is Galois,
the subgroup $\nrm_{k_v}^{k'_w} V(C_{k'_w})$ of $V(C_{k_v})$ 
is the same for all $w | v$.
Thus we may replace the right hand side by 
$\oplus_v [\oplus_{w|v} V(C_{k_v}) 
  /\nrm_{k_v}^{k'_{w}} V(C_{k'_w})]_G$.
(Here $(-)_G$ denotes the coinvariants by $G$,
where $G$ acts on the set $\{ w|v \}$ transitively.)

\noindent
(iii)
The surjectivity of $(**)$ is an easy consequence of
the Moore-Bloch exact sequence
(\cite{bloch} Theorem 4.2, \cite{ks} Proposition 5),
therefore, by \cite{somekawa} Theorem 4.1,
it remains true when one replaces $V(C)$ by $K(k; G, \mathbb{G}_m)$ 
with any semi-abelian group $G$ over $k$.
However, for the proof of the injectivity
we need some deeper results from Kato's Hasse principle \cite{kato},
which is not known for general $K(k; G, \mathbb{G}_m)$.

\noindent
(iv)
The classical Hasse norm theorem states the
injectivity of
\[ k^* /\nrm_k^{k'} (k'^*) \to 
  \oplus_v k_v^* /\nrm_{k_v}^{k'_{w(v)}} (k'^*_{w(v)})
\]
for any cyclic extension $k'/k$ of global fields.
An analogous statement for $K_2$ is known to be valid
for any finite extension $k'/k$
(cf. Bak-Rehmann \cite{br} and Colliot-Th\'el\`ene \cite{ct}).
See also \cite{koya, ost} for related results.
\end{remark}

\begin{proof}
Let $J$ be the Jacobian variety of $C.$
We write $n=[k':k]$ and $M=J[n] \otimes \mu_n.$
We consider the Hochschild-Serre spectral sequence
$E_2^{i, j} = H^i(k, H^j(C_{\bar{k}}, \mu_n^{\otimes 2}))
\Rightarrow H^{i+j}(C, \mu_n^{\otimes 2})$.
We have $E_2^{i, j}=0$ if $j>2$ because $C$ is a curve.
The edge homomorphisms
$E_2^{i, 0} = H^i(k, \mu_n^{\otimes 2}) \to H^i(C, \mu_n^{\otimes 2})$
and
$H^i(C, \mu_n^{\otimes 2}) \to E_2^{i-2, 2} = H^{i-2}(k, \mu_n)$
are split by a $k$-rational point $P$ on $C.$
(The latter is split by the composition
$H^{i-2}(k, \mu_n) \cong H^i_P(C, \mu_n^{\otimes 2}) \to 
H^i(C, \mu_n^{\otimes 2})$
of the purity isomorphism
and the canonical map in the localization sequence.)
As a consequence, 
we have a direct sum decomposition $(**)$
in the following commutative diagram:
\[ 
\begin{matrix}
&SK_1(C)/n& \overset{(*)}{\cong} &k^*/n& \oplus &V(C)/n&
\\
& \downarrow^{f} & & \downarrow^{\cong} & & \downarrow^{g}&
\\
&H^3(C, \mu_n^{\otimes 2})& 
\overset{(**)}{\cong} &H^1(k, \mu_n)& \oplus &H^2(k, M)&
 \oplus H^3(k, \mu_n^{\otimes 2}).
\end{matrix}
\]
Here the decomposition $(*)$ is given by Theorem \ref{rs}.
Both of the maps $f$ and $g$ are injective
(cf. \cite{kato} Lemma 2.4, \cite{yama} Theorem 6.1).
We set $\tilde{D}_n(C) = \coker(f), D_n(C) = \coker(g)$,
so that we have an isomorphism
\[ \tilde{D}_n(C) \cong D_n(C) \oplus H^3(k, \mu_n^{\otimes 2}). \]

We consider a commutative diagram with exact rows
\[
\begin{matrix}
&( V(C_{k'})/n )_G& \to &H^2(k', M)_G& \to &D_n(C_{k'})_G \to 0 \\
&\downarrow^{\nrm_k^{k'}} &  
&\downarrow^{\cre_k^{k'}} &  
&\downarrow^{\phi_k^{k'}}   \\
0 \to & V(C)/n& \to &H^2(k, M)& \to &D_n(C) \to 0,
\end{matrix}.
\]
Similarly, we have a commutative diagram for any $v$
\[
\begin{matrix}
&[\oplus_{w|v} V(C_{k'_w})/n ]_G& 
\to &[\oplus_{w|v} H^2(k'_w, M)]_G& 
\to &[\oplus_{w|v} D_n(C_{k'_w})]_G \to 0 \\
&\downarrow^{\nrm_{k_v}^{k'}} &  
&\downarrow^{\cre_{k_v}^{k'}} &  
&\downarrow^{\phi_{k_v}^{k'}}   \\
0 \to & V(C_{k_v})/n& \to &H^2(k_v, M)& \to &D_n(C_{k_v}) \to 0.
\end{matrix}
\]
Consequently we have a commutative diagram with exact rows
\[
\begin{matrix}
&\ker(\phi_k^{k'}) \to 
&\coker(\nrm_k^{k'}) \to 
&\coker(\cre_k^{k'})
\\
&{}^{a} \downarrow  
&{}^{ } \downarrow  
&{}^{b} \downarrow 
\\
\oplus_v \ker(\cre_{k_v}^{k'}) \to 
&\oplus_v \ker(\phi_{k_v}^{k'}) \to 
&\oplus_v \coker(\nrm_{k_v}^{k'}) \to
&\oplus_v \coker(\cre_{k_v}^{k'}).
\end{matrix}
\]
Since $\coker(\nrm_k^{k'}) = V(C)/\nrm_k^{k'} V(C_{k'})$
(and similarly for $k'_w/k_v$),
the proposition follows from the following lemma.
\end{proof}

\begin{lemma}
(i) $\cre_{k_v}^{k'}$ is an isomorphism
for any $v$.

(ii) $a$ is an isomorphism.

(iii) $b$ is an isomorphism.
\end{lemma}
\begin{proof}
(i) The map $\cre_{k_v}^{k'}$ is Tate dual to
the restriction map
\[ H^0(k_v, M^*) \to (\oplus_{w|v} H^0(k'_w, M^*) )^G  
\qquad (M^* = \homo(M, \mu_n))
\]
which is an isomorphism.

\noindent
(ii) This is a direct consequence of a deep result 
due to Kato \cite{kato} p. 161-163, which says:
\[ \tilde{D}_n(C) \cong \oplus_v \tilde{D}_n(C_{k_v}), \qquad 
\tilde{D}_n(C_{k'}) \cong \oplus_v \oplus_{w|v} \tilde{D}_n(C_{k'_w}). 
\]

\noindent
(iii) This can be shown by the same argument as
Colliot-Th\'el\`ene \cite{ct} Lemma 2 (c).
We include the proof here for the completeness.
Let $M_1 = Ind_k^{k'} M$ be the induced Galois module over $k$
so that $H^*(k, M_1) \cong H^*(k', M)$ by Shapiro's lemma.
We have a surjective map $M_1 \to M$ which induces
$\cre_k^{k'}: H^2(k', M) \to H^2(k, M)$
under the above identification,
and we write $M_2$ for its kernel.
Then we have a commutative diagram with exact rows:
\[
\begin{matrix}
0 \to &\coker(\cre_k^{k'})& \to &H^3(k, M_2)& \to &H^3(k, M_1)& 
\\
& \downarrow &  & \downarrow & &\downarrow & \\
0 \to &
\oplus_v \coker(\cre_{k_v}^{k'_{w(v)}})& 
\to &
\oplus_v H^3(k, M_2)& \to &
\oplus_v H^3(k, M_1)& \\
\end{matrix}
\]
where $w(v)$ is an arbitrary place of $k'$ above $v.$
The theorem of Poitou-Tate (see \cite{serre} \S 6.3 Th\'eor\`eme B)
shows that
the middle and right vertical arrows are isomorphisms,
hence so is the left one.
Since there is an isomorphism
$\coker(\cre_{k_v}^{k'_{w(v)}}) \cong \coker(\cre_{k_v}^{k'})$
for each $v$, we are done.
\end{proof}

\subsection{Example}
\begin{lemma}
Let $E$ be an elliptic curve over a number field $k$.
Let $k'/k$ be a finite Galois extension.
Let $N_1$ be the number of finite places $v$ of $k$ such that
$E_{k_v}$ is a non-split Tate curve $k_v$
and such that 
$E_{k'_w}$ is a Tate curve over $k'_w$
for a place $w$ of $k'$ over $v$
(cf. \S 3.3).
Let $N_2$ be the number of real places $v$ of $k$ such that
$E(k_v)$ has two components 
and such that 
an extension $w$ of $v$ to $k'$ is a complex place.
Then we have
\[  V(E)/\nrm_k^{k'} V(E_{k'}) \cong (\Z/2\Z)^{\oplus (N_1 + N_2)}. \]
\end{lemma}
\begin{proof}
By Proposition \ref{hasse}, we are reduced to compute
$V(E_{k_v})/\nrm_{k_v}^{k'_{w(v)}} V(E_{k'_{w(v)}})$ 
for each place $v$ of $k.$
For a finite place $v,$
this is done by Proposition \ref{normmap} (i),
Lemma \ref{counterlemma}, and Remark \ref{counterremark} (ii).
For an infinite place, this follows from 
a result of Pedrini and Weibel
(\cite{pw1} Lemma 5.6 and Proposition 5.7).
\end{proof}

By the same argument as Corollary \ref{divlemma},
we obtain the following.

\begin{corollary}\label{lastentry}
Let $C$ be a smooth projective geometrically connected curve
over a number field $k$ such that $C(k) \not= \phi.$
Let $E$ be an elliptic curve over $k.$
Then $\tilde{V}(C \times E)/\tilde{V}(C \times E)_{\divi}$
is annihilated by two.
\end{corollary}

\vspace{3mm}
\noindent
{\it Acknowledgement.}
The author would like to thank Kanetomo Sato,
Michael Spiess, and the referee 
for their comments for the earlier version.


\end{document}